\documentclass[]{{birkart}}

\usepackage{graphicx}
\usepackage{amsfonts}
\usepackage{amsmath,amscd}
\usepackage{amsthm}

{ \theoremstyle{plain}
\newtheorem{theorem}{Theorem}
\newtheorem{lemma}{Lemma}

}

{ \theoremstyle{definition}
\newtheorem{definition}{{\bf Definition}}}

{ \theoremstyle{definition}

}

{ \theoremstyle{definition}
\newtheorem{remark}{Remark}
}

\def\R#1{\mathbb{R}^{#1}}
\def\scal#1#2{\langle #1, #2\rangle}
\def\M{\mathcal{M}}

\def\frac#1#2{{#1\over #2}}

\def\uu{\mathbf{u}}

\DeclareMathOperator{\modd}{\mathrm{mod}}

\begin{document}

\title{Life-time of minimal tubes and coefficients
of univalent functions in a circular ring
}

\author{Vladimir \ G.\ Tkachev}
\thanks{This paper was supported in part by Russian Fundamental
Researches Fund, project 93-011-176\\
\textbf{Published in}:
\textit{Bull. Soc. Sci. Lett. {\L}\'{o}d\'{z}. S\'{e}r. Rech. D\'{e}form.} \textbf{20} (1995), 19--26.}

\subjclass{Primary 53A10; Secondary 53A35}

\begin{abstract}
We obtain various estimates of the life-time of two-dimensional minimal tubes in
$\R{3}$ by potential theory methods. 
\end{abstract}

\maketitle

\section{Introduction.}

Let $x=(x_1, x_2,\ldots , x_n, x_{n+1})$ be a point in Euclidean space
$\R{n+1}$ with the {\it time} axis $Ox_{n+1}$ and
$M$ be a $p$-dimensional Riemannian manifold, $2\leq p\leq n$.

\begin{definition}
We say that a surface $\M = (M,\uu )$ given by $C^2$-immersion
$\uu :M\to \R{n+1}$ is a {\it tube} with the projection
interval
$\tau (\M )\subset Ox_{n+1}$, if
(i) for any $\tau\in\tau(\M )$ the
sections $\Sigma_\tau = f(\M)\cap \Pi_\tau$ by hyperplanes
$\Pi_\tau =\{ x\in \R{n+1}_1 : x_{n+1}=\tau\}$ are not empty
compact sets;
(ii) for $\tau', \tau'' \in \tau(\M)$ any part of $\M$ situated
between two different $\Pi_{\tau'}$ and $\Pi_{\tau''}$ is a
compact set.  
\end{definition}

\begin{definition}
A surface $\M$ is called {\it minimal}\/  if the mean curvature of
$\M$ vanishes everywhere.
\end{definition}

It is the well known fact (see \cite{b5}, p.331) that the minimality
condition of $\M$ is equivalent to that all coordination
functions of the immersion $\uu $ are harmonic.  For this
reason, the two-dimensional minimal tubes can be considered as
direct analog of the closed relative string conception in the
modern nuclear physics (cf. \cite{b2}). This approach was proposed by
V.M.Miklyukov and the author in \cite{b7} for an arbitrary dimension
$p$.

From this point of view many intrinsic geometric
invariants of $\M$ have  the natural physical meaning. Namely,
the length of the projection interval $|\tau(\M)|$ can be
interpreted as a {\it life-time}\/ of the tube $\M$.

To introduce the following important characteristic we denote by
$\nu$ the unit normal to $\Sigma_\tau$ with respect to $\M$
which is co-directed with the time-axis $Ox_{n+1}$. Then by virtue
of the harmonicity of the coordinate functions $u_k(m)=x_k\circ
\uu (m)$, $1\leq k\leq n+1$, the flow integrals
$$
J_k=\int_{\Sigma_\tau}\scal{\nabla u_k}{\nu}\;d\Sigma
$$
are independent of $\tau\in\tau(\M)$. Here $d\Sigma$ is
the 1-Hausdorff measure along $\Sigma_\tau$.

\begin{definition}
We call $Q(\M)=(J_1,J_2,\ldots,J_{n+1})\in\R{n+1}$ the \textit{full
flow-vector} of $\M$.
\end{definition}

We notice  the positiveness of $J_{n+1}$ as a consequence of
the choice of $\nu$ direction. Moreover, $Q(\M)$ is an
1-homogeneous functional of $\M$ under the homotheties group
action in $\R{n+1}$. Let us denote by $\alpha(\M)$ the angle
between $Q(\M)$ and the time-axis $Ox_{n+1}$.

In this paper we are interested in the following question:
What sufficient conditions yield the finiteness of the
time-life of a two-dimensional minimal tube?
As it shown in the series of papers \cite{b6}--\cite{b8},
in the case $p\geq 3$ this quantity is always finite and the
following estimation holds
$$
|\tau(\M)|\leq \varrho(\M)c_p,
$$
where $c_p$ depends only on $p$, and $\varrho(\M)$ is the
smallest diameter of sections $\Sigma_\tau$. The last relationship
is sharp and the equality occurs if and only if $\M$ is a
minimal surface of revolution.

A special feature of the two-dimensional case is that there
exist tubes with finite as well as infinite values of the
life-time.  Indeed, a family 
of slanted minimal surfaces with circular
cross-sections $\Sigma_\tau$ was discovered by B.~Riemann \cite{b10}. Some other
recent examples can also be found in \cite{b4}.

In this paper we prove

\begin{theorem}\label{the1}
Let $\M$, ${\rm dim} \M=2$ be a minimal two-connected tube
with univalent Gaussian mapping. If the angle $\alpha(\M)$ is
different from zero, then the life-time $|\tau(\M)|$ of $\M$ is
finite and
$$
\tau(\M)\leq \frac{\pi \|Q\|\cos \alpha(\M)}{\ln \tan
(\frac{\pi}{4}+\frac{\alpha}{2})}.
$$
\end{theorem}

Let us denote by $a_0[f]$ the central
coefficient of the Laurent decomposition of an holomorphic function
$f(z)$ in an annulus $K_R=\{z: 1/R<|z|<R\}$, i.e.
$$
a_0[f]\equiv\int_{C_1}\frac{g(\zeta)\, d \zeta}{\zeta}.
$$
where $C_1$ is the unite circle $\{z\in {\Bbb C}: |z|=1\}$.
The following auxiliary assertion is a key ingredient in the
proof of Theorem~1.

\begin{theorem}\label{the2}
Let $g(z)$ be a univalent holomorphic function defined in the
annulus $K_R$ omitting zero.
Assume that
$$
a_0[g]=\lambda, \qquad a_0[1/g]=-\lambda,
\leqno{(1)}
$$
for some real positive $\lambda$.
Then
$$
\ln R\leq\ln R_0(\lambda)=\frac{\pi^2}{\ln(\lambda+\sqrt{1+\lambda^2})}
\leqno{(2)}
$$
\end{theorem}

\begin{remark} We note that estimate (2) has well
asymptotic behaviour for $R\to\infty$ as shows Riemannian
example mentioned above. But we cannot now present the sharp
value of $R_0(\lambda)$. Nevertheless, it seemed us very
probably that the following conjecture is true.

\end{remark}

\begin{remark} The best upper bound in the left side of (2)
is achieved for holomorphic function
$g_0(z)$ which provides a conformal map of the annulus $K_R$ onto the plain ${\Bbb C}$ with
two slits: $(-1/\alpha;0)$ and $(\alpha;+\infty)$, for the
suitable choice of parameter $\alpha$.

\end{remark}

\noindent
\textbf{Acknowledgement. }
I wish to thank  V.M.Miklyukov for many useful discussions concerning the topic of this paper.

\section{Proof of Theorem~\ref{the2}}

Let $\Gamma=\{C_\rho: 1/R<\rho<R\}$ be a family of all
concentric circles $C_\rho=\{z:|z|=\rho\}$ in the annulus $K_R$.
It follows easily from the non-vanishing property of $g(z)$ that
the loop $C_1$ in the integrals (1) may be replaced by an
arbitrary circle $C_\rho\in\Gamma$. It follows from the mean
value theorem and (1) that for every $\rho\in(1/R;R)$ there
exist $t_1$ and $t_2$ such that
$$
{\rm Re}\; g(\rho
e^{it_1})=\lambda
\quad\hbox{and}\quad
{\rm Re}\; \frac{1}{g(\rho e^{it_2})}=-\lambda.
\leqno{(3)}
$$

Let $\gamma_\rho=g(C_\rho)$. Then by virtue of the univalence
of $g(z)$, the curve $\gamma_\rho$ is the simple Jordan one. Let
$g(\rho e^{it})=x(t)+iy(t)$ be the representation of
$\gamma_\rho$. Then we obtain from (3)
$$
x(t_1)=\lambda;\qquad x^2(t_2)+y^2(t_2)+\frac{1}{\lambda}x(t_2)=0.
$$

The last relations have the helpful geometric interpretation:

\medskip
\noindent
($\star$) The curve $\gamma_\rho$ has a non-empty intersection with the vertical rightline
$L_1=\{z: {\rm Re}z=\lambda\}$ and the circle $L_2=\{z:
|z+1/2\lambda|=1/2\lambda\}$.
\medskip

We shall make use the technique  from the potential
theory (the length-are method). Recall the exact definition.
Let $E$ be a family of locally rectifiable curves $\gamma$ and
$\varphi(z)\geq 0$ be a Baire function with the property
$$
\int_\gamma \varphi(z)\,|dz|\geq 1,
$$
for every  $\gamma\in E$. The infimum
$$
\modd  \;E=\inf \int \varphi^2(z)\, dx\,dy
$$
over all such $\varphi(z)$ is called a {\it conformal module }
of the family $E$.

Then it is known (see \cite{b1}) that $\modd \, E$ is the conformal
invariant. As a consequence we obtain in our  situation
$$
\modd \,\Gamma=\modd \, \Gamma_1,
\leqno{(4)}
$$
where $\Gamma_1=\{\gamma_\rho: 1/R<\rho<R\}$.

Let us denote by $D$ the two-dimensional domain
$$
D=\left\{z: {\rm Re}\,z<\lambda;\;
\left|z+\frac{1}{2\lambda}\right|>
\frac{1}{2\lambda}\right\}.
$$
Using the ($\star$)-property, we can find for every $\rho\in
(1/R;R)$ the continuum $\gamma'_\rho\subset\gamma_\rho$
joining the boundary components of $D$.
Then a family $\Gamma_2$ consisting of all continua
$\gamma'_\rho$ is ``shorter'' than $\Gamma_1$ and
it follows from Theorem 1.2, \cite{b1} that
$$
\modd \, \Gamma_1\leq\modd \, \Gamma_2.
\leqno{(5)}
$$

On the other hand, $\Gamma_2$ is the subfamily of $\Gamma(D)$,
where the last term means the family of {\it all} curves joining
the boundary components of a domain $D$. The monotonicity
property of infimum and Definition 4 lead to the following
inequality
$$
\modd \, \Gamma_2\leq\modd \Gamma(D).
\leqno{(6)}
$$

Now, combining the standard fact
$$
\modd \, \Gamma =\frac{\ln R}{\pi}
\leqno{(7)}
$$
with relations (4), (5) and (6) we arrive at the following
inequality
$$
\frac{\ln R}{\pi}\leq \modd \Gamma(D).
$$

To compute the last module we note that the linear-fractional
function
$$
f(z)=\frac{1}{\lambda^*}\cdot\frac{z+\lambda^*}{1-z\lambda^*}
$$
maps $D$ onto an annulus $K_1=\{w: 1<|w|<1/{\lambda^*}^2\}$,
where $\lambda^*=\sqrt{\lambda^2+1}-\lambda$. Thus, using the
invariance property of conformal module we obtain
$$
\frac{\ln R}{\pi}\leq \modd (D)\equiv \frac{2\pi}
{\ln{(1/\lambda^*}^2)} =
\frac{\pi}{\ln (\lambda+\sqrt{1+\lambda^2})}.
$$
and Theorem 2 is proved.

\section{The Gaussian map of two-dimensional minimal tubes
and their full-flow vector}

In this section we express the full flow-vector of
an arbitrary two-dimensional tube $\M\in\R{n}$ via
Chern-Weierstrass representation for minimal surfaces. Namely,
if $\M$ is a two-connected surface then we can arrange that
$\M$ is conformally equivalent to an annulus $K_R$ for the
appropriate $R>1$. Then there exist the corresponding
parametrization of $\M$ (see \cite{b9}):
$$
\uu (z)={\rm Re}\;\int_{z_0}^{z}F(\zeta)\,d\zeta\;:
K_R\rightarrow \R{n},
$$
where
$$
F(z)=(\varphi_1(\zeta),\ldots, \varphi_n(\zeta))
$$
and $\varphi_i(\zeta)$ are holomorphic functions satisfying the
following conditions
$$
\sum_{i=1}^{n}\varphi_i(\zeta)^2=0;
\leqno{(8)}
$$
and
$$
{\rm Re}\, \int_{|z|=1}F(\zeta)\,d\zeta={\bf 0}.
\leqno{(9)}
$$

\begin{lemma}\label{lem1}
Under the above hypotheses we have
$$
Q(\M)={\rm Im} \int_{|z|=1}F(\zeta)\,d\zeta.
\leqno{(10)}
$$
\end{lemma}

\begin{proof}
It sufficient to show that
$$
J_k\equiv\int_{\Sigma_\tau}\scal{\nabla u_k}{\nu}\;d\Sigma= {\rm
Im}  \int_{|z|=1}\varphi_k(\zeta)\,d\zeta,
\leqno{(11)}
$$
for every $k=1,2,\ldots, n+1$.

To prove (11) we introduce the conjugate to $u_k(z)$ function
$v_k(z)$ by
$$
v_k^*(z)={\rm Im}\int_{z_0}^z \varphi_k(\zeta)\,d\zeta,
$$
We notice that $v_k(z)$ in general is a multivalued function.
On the other hand, the covariant derivative $\nabla v_k$ is well
defined and using the properties of Hodge $\star$~-~operator we
have
$$
\int_{\Sigma_\tau}\scal{\nabla u_k}{\nu}\;d\Sigma=
\int_{\Sigma_\tau}\scal{\star\nabla u_k}{\star\nu}\;d\Sigma=
\int_{\Sigma_\tau}\scal{\nabla v_k}{\star\nu}\;d\Sigma=
$$
$$
=\int_{\Sigma_\tau} d\,v_k=
{\rm Im}\int_{|z|=1}\varphi_k(\zeta)\, d\zeta,
$$
and (11) is proved.
\end{proof}

In our case $n=2$, Chern-Weierstrass representation can be
simplified in the following classic way. Namely, there exist
a holomorphic function $f(z)$ and a meromorphic function $g(z)$
which are well defined in the annulus $K_R$ and such that
$$
F(z)=\left((1-g^2)f;\, i(1+g^2)f;\, 2gf\right).
\leqno{(12)}
$$
Moreover, poles of $g(z)$ coincide with zeros of $f(z)$ and
the order of a pole of $g(z)$ is precisely the order of the
corresponding zero of $f(z)$.
We emphasize that $g(z)$ is a composition of the stereographic
projection and Gaussian map of $\M$.

\begin{lemma}\label{lem2}
In our assumptions
$$
2fg\equiv\frac{\scal{Q(\M)}{e_3}}{2\pi z},
\leqno{(13)}
$$
and $g(z)$ omits the zero and infinity values.
\end{lemma}

\begin{proof}
We use the method proposed by M.~Schiffman in \cite{b11}.
We recall that the coordinate function $u_3(z)$ is harmonic in
the annulus $K_R$ and by virtue of Definition 1,
$$
\lim_{z\to 1/R }u_3(z)=\tau_1,
\qquad
\lim_{z\to R }u_3(z)=\tau_2,
\leqno{(14)}
$$
where $\tau(\M)=(\tau_1;\tau_2)$ is the projection of the tube
$\M$ onto $x_3$-axis.

We consider an auxiliary harmonic function
$$
h(z)=\tau_1+\frac{\tau_2-\tau_1}{2\ln R}\ln |z|.
$$
It is easily seen that $h(z)$ satisfies (14). Thus
$h_1(z)=u_3(z)-h(z)$ is harmonic in the annulus and
$$
\lim_{z\to\partial K_R}h_1(z)=0.
$$
Then the maximum principle implies that $h_1(z)\equiv 0$
everywhere in $K_R$ and hence
$$
u_3(z)\equiv \tau_1+\frac{\tau_2-\tau_1}{2\ln R}\ln |z|.
\leqno{(15)}
$$

In particular, it follows from (15) that
$$
d u_3(z)\equiv
\frac{\tau_2-\tau_1}{\ln R}\cdot\frac{z}{|z|^2}
$$
doesn't vanish in $K_R$. We have, as a consequence, the normal
$n(z)$ to $\M$ isn't parallel to $e_3$ at any point. Taking into
account the above remark about the geometrical sense of $g(z)$ we
obtain that $g(z):K_R\to {\Bbb C}-\{0;\infty\}$.

By comparing of (15) and (12) we deduce that
$$
2g(z)f(z)=\frac{\tau_2-\tau_1}{2\ln R}\cdot \frac{dz}{z}.
\leqno{(16)}
$$
In order to eliminate $\ln R$ from the latter equality we substitute (16)
into (12), and after using (10) we obtain
$$
\ln R=\frac{\pi(\tau_2-\tau_1)}{J_3}.
\leqno{(17)}
$$

On substituting of the found relationship into (16) we arrive at
the conclusion of the lemma.
\end{proof}

\section{Proof of Theorem \ref{the1}}

Let us denote $w=(J_1+iJ_2)/J_3$. Combining Lemma 2, (12) and
(9) we obtain
$$
\int_{C_1}\frac{1-g^2(\zeta)}{2g(\zeta)}\frac{d\zeta}{\zeta}=
2\pi w_1i,
$$

$$
\int_{C_1}\frac{1+g^2(\zeta)}{2g(\zeta)}\frac{d\zeta}{\zeta}=
2\pi w_2.
$$
Simplifying the last expressions and denoting $w=|w|\cdot
e^{i\theta}$, $g_1(z)=-e^{-i\theta}g(z)$ give the following
system
$$
\frac{1}{2\pi}\int_{C_1}\frac{g_1(\zeta)d\zeta}{\zeta}=|w|,
$$

$$
\frac{1}{2\pi}\int_{C_1}\frac{d\zeta}{g_1(\zeta)\zeta}=-|w|.
$$
Applying Theorem 2 we arrive at the inequality
$$
\ln R\leq \frac{\pi^2}{|w|+\sqrt{1+|w|^2}}
$$
where $|w|\equiv |J_1+i J_2|/|J_3|=\tan \alpha(\M)$. Using (17)
we obtain the required estimate and the theorem is proved.


\end{document}